\documentclass[11pt]{article}
\usepackage{amsmath, amsthm}
\usepackage{amscd}
\usepackage{amsfonts}
\usepackage{amssymb}
\usepackage{float,comment}
\usepackage{mathrsfs}
\usepackage{tikz}
\usetikzlibrary{patterns}
\usepackage{color}
\definecolor{red}{rgb}{.7,0,0}
\definecolor{blue}{rgb}{0,0,1}

\usepackage[labelfont=bf,labelsep=space]{caption}

\newtheorem{theorem}{Theorem}[section]
\newtheorem{proposition}[theorem]{Proposition}
\newtheorem{lemma}[theorem]{Lemma}

\theoremstyle{definition}

\def\proof{{\noindent\sc Proof. \quad}}

\newcommand\proofof[1]{{\noindent\sc Proof of #1. \quad}}
\def\eproof{{\mbox{}\hfill\qed}\medskip}
\begin{document}

\makeatletter

\renewcommand{\bar}{\overline}

%Other math symbols
\newcommand{\x}{\times}
\newcommand{\<}{\langle}
\renewcommand{\>}{\rangle}
\newcommand{\into}{\hookrightarrow}

%Greek letters
\renewcommand{\a}{\alpha}
\renewcommand{\b}{\beta}
\renewcommand{\d}{\delta}
\newcommand{\D}{\Delta}
\newcommand{\e}{\varepsilon}
\newcommand{\g}{\gamma}
\newcommand{\G}{\Gamma}
\renewcommand{\l}{\lambda}
\renewcommand{\L}{\Lambda}
\newcommand{\n}{\nabla}
\newcommand{\var}{\varphi}
\newcommand{\s}{\sigma}
\newcommand{\Sig}{\Sigma}
\renewcommand{\t}{\theta}
\renewcommand{\O}{\Omega}
\renewcommand{\o}{\omega}
\newcommand{\z}{\zeta}
\newcommand{\balpha}{\boldsymbol \alpha}
\newcommand{\ab}{\alpha_{\bullet}}
\def\ba{\ol{\a}}
\def\bb{\ol{\b}}
\def\bg{\ol{\gamma}}
\def\bt{\ol{\tau}}
\def\mun{\mu_{\mathrm{norm}}}
\def\sep{\mathsf{sep}}
\def\vol{\mathsf{vol}}
\def\scap{\mathsf{cap}}
\def\Tub{\mathrm{Tub}}
\def\hmg{\mathrm{proj}}
\def\uno{\mathsf{1\hspace*{-2pt}l}}
\def\sfv{\mathsf{v}}

%Other macros
\newcommand{\p}{\partial}
\renewcommand{\hat}{\widehat}
\renewcommand{\bar}{\overline}
\renewcommand{\tilde}{\widetilde}
\newcommand{\diff}{\mathrm{d}}
\newcommand{\erf}{\mathrm{erf}}
%%%%%%%%
% Some fonts
%%%%%%%

%\newcommand{\fiverm}{\tiny\rm}
\font\eightrm=cmr8
\font\ninerm=cmr9

%%%%%
% The next macros define fonts for reals, rationals, complex,
% integers and natural numbers by \R,\Q,\C,\Z and \N respectively.
% Also, \Ri gives \R to the infinity.
%%%%%

\def\N{\mathbb{N}}
\def\Z{\mathbb{Z}}
\def\R{\mathbb{R}}
\def\Q{\mathbb{Q}}
\def\C{\mathbb{C}}
\def\F{\mathbb{F}}
\def\IS{\mathbb{S}}
\def\proj{\mathbb{P}}
\def \Ri{\R^\infty}
\def \Zi{\Z^\infty}
\def \ZRi{\Z^\infty\x\R^\infty}
\def \SZi{\Z^\infty\x S(\R^\infty)}

%%%%%%%
% For writting algorithms
%%%%%%%

\newcommand{\algoritmo}{\begin{minipage}{0.87\hsize}\linea}
\newcommand{\falgoritmo}{\linea\end{minipage}\bigskip}
\newcommand{\codigo}{\begin{minipage}{0.87\hsize}}
\newcommand{\fcodigo}{\end{minipage}\bigskip}
\newcommand{\linea}{\vspace*{-5pt}\hrule\vspace*{5pt}}
\newtheorem{algorithm}{Algorithm}
\def\espacio{\hspace*{1cm}}
\def\eespacio{\hspace*{1.5cm}}
\def\eeespacio{\hspace*{2cm}}
\newcommand{\inputalg}[1]{\linea\bf Input:\quad\rm #1\vspace*{3pt}}
\newcommand{\specalg}[1]{\bf Preconditions:\quad\rm #1}
\newcommand{\Output}[1]{\linea\bf Output:\quad\rm #1\vspace*{2pt}}
\newcommand{\postcond}[1]{\bf Postconditions:\quad\rm #1\vspace*{3pt}}
\newcommand{\bodyalg}[1]{\linea\tt #1\vspace*{3pt}} \newcommand{\bodycode}[1]{\tt #1\vspace*{3pt}}
\newcommand{\sssection}[1]{\ \\ \noindent{\bf #1.}\quad}

%%%%%%%
% Some mathematical operators and constants
%%%%%%%

\def\sg{{\rm sign}\,}
\def\sig{\overline{\rm sign}\,}
\def\rank{{\rm rank}\,}
\def\Im{{\rm Im}}
\def\dist{{\rm dist}}
\def\degree{{\rm degree}\,}
\def\grad{{\rm grad}\,}
\def\size{{\rm size}}
\def\supp{{\rm supp}}
\def\Id{{\rm Id}}
\def\spann{\mathop{\rm span}}
\def\fl{\mathop{\tt fl}}
\def\op{\mathop{\tt op}}
\def\cost{{\rm cost}}
\def\emac{\varepsilon_{\mathsf{mach}}}
\def\cond{\mathop{\mathsf{cond}}}
\def\card{{\rm card}}
\def\mod{{\rm mod}\;}
\def\trace{{\rm trace}}
\def\Prob{\mathop{\rm Prob}}
\def\length{{\rm length}}
\def\Diag{{\rm Diag}}
\def\diam{{\rm diam}}
\def\Tr{{\mathsf{Tr}}}
\def\arg{{\rm arg}\,}
\def\ot{\leftarrow}
\def\transp{^{\rm t}}
\def\bL{{\bf L}}
\def\bR{{\bf R}}
\def\bC{\mathbf{C}}
\def\bfE{{\bf E}}
\def\bd{{\bf d}}
\def\adj{{\rm adj}\:}
\def\cruz{\raise0.7pt\hbox{$\scriptstyle\times$}}
\def\todif{\stackrel{\scriptscriptstyle\not =}{\to}}
\def\nedif{\stackrel{\scriptscriptstyle\not=}{\nearrow}}
\def\sedif{\stackrel{\scriptscriptstyle\not=}{\searrow}}
\def\nodown{\downarrow\mkern-15.3mu{\raise1.3pt\hbox{$\scriptstyle\times$}}}
\def\noto{\to\mkern-20mu\cruz}
\def\notox{\to\mkern-23mu\cruz}
\def\alg{{\hbox{\scriptsize\rm alg}}}
\def\ttl{{\tt l}}
\def\tto{{\tt o}}
\def\ttml{\bar{\tt l}}
\def\ttmo{\bar{\tt o}}
\def\ri{R^\infty}
\def\Oh{{\cal O}}
\def\parts{{\cal P}}
\def\coeff{{\hbox{\rm coeff}}}
\def\Gr{{\hbox{\rm Gr}}}
\def\Error{{\hbox{\tt Error}}\,}
\def\appleq{\hbox{\lower3.5pt\hbox{$\;\:\stackrel{\textstyle<}{\sim}\;\:$}}}

% Para usar en dibujos con Pictex %%
\def\bolita{\scriptscriptstyle\bullet}

%%%%%%%
% Some abreviations for the bib database
%%%%%%%

\def\JACM{Journal of the ACM}
\def\CACM{Communications of the ACM}
\def\ICALP{International Colloquium on Automata, Languages
            and Programming}
\def\STOC{annual ACM Symp. on the Theory
          of Computing}
\def\FOCS{annual IEEE Symp. on Foundations of Computer Science}
\def\SIAM{SIAM J. Comp.}
\def\SIOPT{SIAM J. Optim.}
\def\BSMF{Bulletin de la Soci\'et\'e Ma\-th\'e\-ma\-tique de France}
\def\CRAS{C. R. Acad. Sci. Paris}
\def\IPL{Information Processing Letters}
\def\TCS{Theoret. Comp. Sci.}
\def\BAMS{Bulletin of the Amer. Math. Soc.}
\def\TAMS{Transactions of the Amer. Math. Soc.}
\def\PAMS{Proceedings of the Amer. Math. Soc.}
\def\JAMS{Journal of the Amer. Math. Soc.}
\def\LNM{Lect. Notes in Math.}
\def\LNCS{Lect. Notes in Comp. Sci.}
\def\JSL{Journal for Symbolic Logic}
\def\JSC{Journal of Symbolic Computation}
\def\JCSS{J. Comput. System Sci.}
\def\JoC{J. Complexity}
\def\MP{Math. Program.}
\sloppy

\bibliographystyle{plain}

%new macros
%%%%%%%%%%%%%%%%%%%

\def\CPRi{{\rm \#P}_{\kern-2pt\R}}

\newcommand{\bfb}{{\boldsymbol{b}}}
\newcommand{\bfd}{{\boldsymbol{d}}}
\newcommand{\bff}{{\boldsymbol{f}}}
\newcommand{\bfx}{{\boldsymbol{x}}}
\newcommand{\bfa}{{\boldsymbol{a}}}
\newcommand{\ol}[1]{\overline{#1}}
\def\msC{\mathop{\mathscr C}}
\def\mcC{{\mathcal C}}
\def\mcD{{\mathcal D}}
\def\mcN{{\mathcal N}}
\def\mcG{{\mathcal G}}
\def\mcU{{\mathcal U}}
\def\mcM{{\mathcal M}}
\def\mcZ{{\mathcal Z}}
\def\mcB{{\mathcal B}}
\def\mcX{{\mathcal X}}
\def\PP{{\mathscr P}}
\def\AA{{\mathscr A}}
\def\scC{{\mathscr C}}
\def\scO{{\mathscr O}}
\def\scE{{\mathscr E}}
\def\sG{{\mathscr G}}
\def\mZ{{\mathcal Z}}
\def\mI{{\mathcal I}}
\def\sH{{\mathscr H}}
\def\sF{{\mathscr F}}
\def\bD{{\mathbf D}}
\def\nrr{{\#_{\R}}}
\def\Sd{{\Sigma_d}}
\def\oPp{\overline{P'_a}}
\def\oP{\overline{P_a}}
\def\oDH{\overline{DH_a^\dagger}}
\def\oH{\overline{H_a}}
\def\ua{\overline{\a}}
\def\Hd{\HH_{{\boldsymbol{d}}}}
\def\Hdm{\Hd[m]}
\def\Lg{{\rm Lg}}
\def\sfC{{\mathsf{C}}}
\def\Rk{{\operatorname{rank}}}
%%%%%

\def\bx{{\bf x}}
\def\ii{{\'{\i}}}

\def\P{\mathbb P}
\def\Exp{\mathop{\mathbb E}}

\newcommand{\binomial}[2]{\ensuremath{{\left(
\begin{array}{c} #1 \\ #2 \end{array} \right)}}}

\newcommand{\HH}{\ensuremath{\mathcal H}}
\newcommand{\diag}{\mathbf{diag}}
\newcommand{\CH}{\mathsf{CH}}
\newcommand{\Cone}{\mathsf{Cone}}
\newcommand{\SCH}{\mathsf{SCH}}

%%%classes

\newcommand{\NPR}{\mathsf{NP}_{\R}}
\newcommand{\NPc}{\mathsf{NP}}
\newcommand{\sfH}{\mathsf{H}}
\newcommand{\SPR}{\#\mathsf{P}_{\R}}
\newcommand{\PSPACE}{\mathsf{PSPACE}}
\newcounter{line}
%\newcounter{algorithm}

\newcommand{\macheps}{\varepsilon_{\mathrm{m}}}
\newcommand{\sgn}{\mathrm{sgn}}

\begin{title}
{\LARGE {\bf On local analysis}}
\end{title}
\author{Felipe Cucker\thanks{Partially
supported by a GRF grant from the Research Grants
Council of the Hong Kong SAR (project number CityU 11302418).}
\\
Dept. of Mathematics\\
City University of Hong Kong\\
{\tt macucker@cityu.edu.hk}
\and
Teresa Krick\thanks{Corresponding author. Partially supported by grant  CONICET-PIP2014-2016-112
20130100073CO.}\\
Departamento de Matem\'atica \& IMAS\\
Univ. de Buenos Aires \&\ CONICET\\
ARGENTINA\\
{\tt krick@dm.uba.ar}
}

\date{}
\makeatletter
\maketitle
\makeatother

\begin{quote}
{\small
{\bf Abstract.}
We extend to Gaussian distributions a result providing smoothed analysis estimates
for condition numbers given as relativized distances to ill-posedness. We also
introduce a notion of {\em local analysis} meant to capture the behavior of these
condition numbers around a point.
}\end{quote}

\section{Introduction}

In the 1990s D.~Spielman and S.H.~Teng introduced the notion of {\em smoothed analysis},
in an attempt to give a more realistic analysis of the practical performance of an algorithm
than those obtained through the use of worst-case or average-case analyses.
In a nutshell, this new paradigm in probabilistic analysis interpolates between
worst-case and average-case by considering the worst-case (over the data)
of the average value (over possible random perturbations) of the analyzed
quantity. See, for instance, \cite{SpiTen2009} for an overview.

An example of this analysis to the quantity $\ln\kappa(A)$, where $A$
is a square matrix and $\kappa(A):=\|A\|\,\|A^{-1}\|$,
was provided by M.~Wschebor in~\cite{Wsch:04}. Wschebor showed that
\begin{equation}\label{eq:Mario}
  \max_{\bar{A}\in\IS(\R^{n\times n})}
  \Exp_{A\sim N(\bar{A},\sigma^2\Id)}\ln\kappa(A) \leq
  \ln \Big(\frac{n}{\min\{\sigma,1\}}\Big) +\Oh(1),
\end{equation}
where here, and in what follows, $x\sim N(\bar{x},\sigma^2\Id)$ indicates that $x$
is drawn from an isotropic  Gaussian distribution centered at $\bar{x}$ with
covariance matrix $\sigma^2\Id$. The behavior of the bound $H(n,\sigma)$ in the
right-hand side of~\eqref{eq:Mario} shows two expected properties of a
smoothed analysis:
\begin{description}
\item[(SA1)]  When $\sigma\to0$, $H(n,\sigma)$  tends to its worst-case value
(there are no random perturbations of the input in this case).
\item[(SA2)] When $\sigma\to\infty$, $H(n,\sigma)$   tends to the average value of the
  analyzed quantity (the random perturbation is over all the input data in this case).
 \end{description}
Indeed, the convergence of $H(n,\sigma)$  to
infinity when $\sigma\to0$ is clear, and with it (SA1).
And a result of A.~Edelman~\cite{Edelman88}
proves that $\Exp_{A\sim N(0,\sigma^2\Id)}\ln\kappa(A)=\ln n +\Oh(1)$,
thus showing (SA2).

The main  agenda of this paper is to
introduce  the notion of {\em local analysis}, which aims to study locally at a
base point $\bar{x}$ the average value over possible random perturbations of the
analyzed quantity, without taking then the worst-case over all input data.
The benefit of such analysis is that it
provides information depending directly on the base point
instead of assuming a worst-case, as in the smoothed analysis.

We illustrate this notion by developing it for a {\em conic condition
number}. This is a condition number satisfying a
Condition Number Theorem. We next describe more precisely this notion
and its context.

In~1936 Eckart and Young~\cite{EckYou} proved that for a square matrix $A$, $\kappa(A)=\|A\|/d(A,\Sigma)$
where $\Sigma$ is the set of non-invertible matrices and $d$ denotes distance.
This result came to be known as the {\em Condition Number Theorem},
even though it was proved more than
ten years before the introduction of condition numbers by
Turing~\cite{Turing48} and von Neumann and Goldstine~\cite{vNGo47}.
In~1987 J.~Demmel observed (and proved) that similar Condition Number Theorems  hold true
for the condition numbers of various problems~\cite{Demmel87}. More precisely,
he showed that these condition numbers were either equal to or closely bounded
by the (normalized) inverse to the distance to ill-posedness. That is,
that for an input data $x$ of the problem at hand, the condition number
of $x$ for that problem is either equal to or closely bounded by
\begin{equation}\label{eq:CCN}
   \msC(x)=\frac{\|x\|}{d(x,\Sigma)},
\end{equation}
where  $\Sigma\ne \{0\}$ is an algebraic cone of {\em ill-posed inputs}.
One year later, Demmel~\cite{Demmel88}
derived general average analysis bounds for those (conic) condition
numbers.
These bounds depend only on the dimension $N+1$ of the ambient
space, the codimension of $\Sigma$, and its degree. He carried
out this idea for the complex case and stated it for the real case
(requiring $\Sigma$ to be complete intersection)
based on an unpublished (and not findable anywhere) result
by Ocneanu. The underlying probability distribution is the
isotropic Gaussian on $\R^{N+1}$ but it is easy to observe that
the bounds hold as well for the uniform distribution on
the unit sphere $\IS^N$  (or, equivalently, on any half-sphere,
due to the equality $\msC(-x)=\msC(x)$).

In~\cite{BuCuLo:07} Demmel's idea was extended to perform a
smoothed analysis of the conic condition number $\msC(x)$
in the case that $\Sigma$ is the zero set of a single real homogeneous
polynomial $F$ in $N+1$ variables. For this analysis one considers
the centers $\bar{x}$ of the distributions in $\IS^N$ (as in~\eqref{eq:Mario})
and there are two natural choices for the distribution itself:
a Gaussian supported in $\R^{N+1}$ or a uniform on a spherical cap in $\IS^N$.
The uniform case is studied in \cite{BuCuLo:07}, where the following
bound is obtained for $\theta\in[0,\pi/2]$:
\begin{equation}\label{eq:smooth-u}
  \max_{\bar{x}\in\IS^N}\Exp_{x\in B_{\IS}(\bar{x}, \theta) } \ln \msC(x)
  \;\le\; \ln \frac{Nd}{\sin\theta} + 2(\ln 2 + 1)
\end{equation}
where $d$ is the degree of $F$  and
$B_{\IS}(\bar{x},\theta)$ is the spherical cap of radius
$\theta$ centered at $\bar{x}$ which we endow with the
uniform distribution. This bound $H(N,d,\theta)$ recovers an
average analysis in the particular case that the spherical cap
is a half-sphere. That is,
\begin{description}
\item[(SA2')]
   $H(N,d,\pi/2)=\ln (Nd)+\Oh(1)$, is the average value of $\ln \msC(x)$ for $x\in \IS^N$, see~\cite{Demmel88}.
\end{description}

A smoothed analysis of the conic condition number $\msC(x)$ in the Gaussian case $N(\bar x, \sigma^2\Id)$  was still lacking, and it is one of the results we present in this paper, since it is strongly linked with our local analysis as we will see below. Theorem  \ref{th:smoothedGaussian} shows that
$$ \max_{\bar x\in \IS^N} \Exp_{x\sim N(\bar x,\sigma^2\Id)}\ln \msC(x) \, \le \, H(N,d,\sigma)$$
where $H(N,d,\sigma)$ is an explicit bound that satisfies {\bf (SA1)} and {\bf (SA2)}. That is
$$\lim_{\sigma\to 0} H(N,d,\sigma)=\infty \quad \mbox{and} \quad
\lim_{\sigma\to \infty} H(N,d,\sigma)=\ln ({Nd}) + \Oh(1).$$

With respect to   local analysis, the gist is to
obtain bounds for the quantities
$$
    \Exp_{x\sim \mcD(\bar{x})}\ln\msC(x)
$$
where $\bar{x}\in\IS^N$ and $\mcD(\bar{x})$ is either the
uniform distribution on the spherical cap $B_{\IS}(\bar{x},\theta)$
or the Gaussian $N(\bar{x},\sigma^2\Id)$.

These bounds
will be expressions $H(N,d,\nu,\msC(\bar{x}))$ where $\nu$
is either $\theta$ or $\sigma$ depending on the underlying
distribution, which should coincide
with smoothed analysis bounds when $\msC(\bar{x})=\infty$.
More precisely, if we denote by
$H_\infty(N,d,\nu)$ the result of replacing $\msC(\bar{x})$ by
$\infty$ in $H(N,d,\nu,\msC(\bar{x}))$ then we want the following:
\begin{description}
\item[(LA0)]
  $H_\infty(N,d,\nu)$ has the same behavior as the smoothed analysis
  bound $H(N,d,\nu)$.
\end{description}
Furthermore, when $\msC(\bar{x})<\infty$ we seek the following
limiting behavior:
\begin{description}
\item[(LA1)]
  $\displaystyle\lim_{\nu\to 0} H(N,d,\nu,\msC(\bar{x})) =
\ln(\msC(\bar{x}))+ \Oh(1)$, the local complexity at $\bar x$.
\item[(LA2)]
  $\displaystyle\lim_{\sigma\to \infty} H(N,d,\sigma,\msC(\bar{x})) =
\ln(Nd)+ \Oh(1)$ in the Gaussian case, the average complexity.
\item[(LA2')]
  $H(N,d,\pi/2,\msC(\bar{x})) =\ln(Nd)+ \Oh(1)$ in the uniform case, the average complexity.
\end{description}

Indeed, we show that this is the case in Theorem~\ref{thm:instance} (uniform
case) and Theorem~\ref{thm:main-local} (Gaussian case).

\bigskip
\noindent {\bf  Acknowledgments.}  We are grateful to Pierre Lairez for many useful discussions. In
particular, for pointing to us an argument in Proposition 4.2.
\section{Notations and preliminaries}

In all what follows we consider the space $\R^{N+1}$
endowed with the standard inner product $\langle~ ,~\rangle $
and its induced norm $\|~\|$. Within this space we have the unit sphere
 $\IS^N=\{x\in \R^{N+1}:\, \|x\|=1\}$, and for $\bar x\in \IS^N$ we denote
by $B(\bar x,r)=\{x\in\R^{N+1}:\, \|x-\bar x\|\leq r\}$ the closed ball centered
at $\bar x\in\R^{N+1}$ with radius $r\ge 0$, and by
$$
 B_{\IS}(\bar x,\theta)=\{x\in\IS^N:\, 0\le \sphericalangle(x,\bar x)\leq \theta\}
 =\{x\in \IS^{N}\,:\,\langle x,\bar x\rangle \ge \cos\theta\}
$$
the spherical cap in $\IS^N$ centered at $\bar x\in \IS^N$ with radius
$0\le \theta\le \pi$, that is the closed ball of radius $\theta$ around $\bar x$
in $\IS^N$ with respect to the Riemannian distance in $\IS^N$.

We will also refer to the sine distance $d_{\sin}$ in
$\R^{N+1}\setminus\{0\}$ given by
$d_{\sin}(x,\bar x):=\sin(\sphericalangle(x,\bar x))$.
Let $B_{\sin}(\overline{x},\rho):=\{x\in \IS^N: \,
d_{\sin}(x,\overline{x})\le \rho\}$ denote the closed ball of
radius $\rho$ with respect to
$d_{\sin}$ around $\bar{x} \in \IS^N$.
This is the union of $B_{\IS}(\bar{x},\theta)$
with $B_{\IS}(-\bar{x},\theta)$ where $\theta\in[0,\pi/2]$ is
such that $\rho=\sin\theta$.

We will denote by $\Oh_N=\vol(\IS^N)$ the volume of $\IS^N$.
We recall (see~\cite[Prop.~2.19(a)]{Condition}) that
\begin{equation}\label{eq:volSN}
 \Oh_N=\frac{2\pi^{\frac{N+1}{2}}}{\Gamma(\frac{N+1}{2})}
\end{equation}
as well as~\cite[Cor.~2.20]{Condition}
\begin{equation}\label{eq:volB}
 \vol(B(0,1))=\frac{\Oh_N}{N+1}
\end{equation}
and, for $x\in\IS^N$ and $\theta\in[0,\frac{\pi}{2}]$,
the bound (see~\cite[Lem.~2.34]{Condition})
\begin{equation}\label{eq:volBS}
  \frac{\Oh_N}{\sqrt{2\pi (N+1)}}(\sin \theta)^N
  \le \vol(B_{\IS}(x,\theta))\leq\frac{\Oh_N}{2} (\sin\theta)^N.
\end{equation}

The main object  in this paper is a
{\em conic condition number} on $\R^{N+1}$, i.e.  a function given by
$$
  \msC:\R^{N+1}\to [1,\infty],\quad
  \msC(x)=\frac{\|x\|}{d(x,\Sigma)},
$$
where  $\Sigma\ne \{0\}$ is the set of {\em ill-posed inputs} in
$\R^{N+1}$, which we assume closed under scalar multiplication.
We note that $\msC(x)\ge 1$ for all $x$ since $0\in \Sigma$.
As $\msC$ is scale invariant
we may restrict to data $x$ lying in $\IS^N$
where $\msC$ can also be expressed as
$$
   \msC(x)=\frac{1}{d_{\sin}(x,\Sigma\cap \IS^N)}.
$$

\section{The uniform case}\label{sec:uniform}

We endow $B_{\sin}(\bar{x},\rho)$ with the uniform probability
measure.  A smoothed analysis for this measure is
given in~\cite[Th. 21.1]{Condition}. Assume that $\Sigma$ is contained
in a real algebraic
hypersurface, given as the zero set of a homogeneous polynomial of
degree~$d$. Then, for all $\theta\in[0,\frac{\pi}2]$ and $\rho:=\sin \theta$,
we have
\begin{equation}\label{eq:smooth}
  \Exp_{x\in B_{\IS}(\bar{x}, \theta) } \ln \msC(x)
  = \Exp_{x\in B_{\sin}(\bar{x}, \rho) } \ln \msC(x)
  \;\le\; \ln \frac{Nd}{\sin\theta} + K
\end{equation}
and
\begin{equation}\label{eq:smooth-pi}
  \Exp_{x\in {\IS}^N } \ln \msC(x)
  \;\le\; \ln (Nd) + K,
\end{equation}
where $K=2(\ln 2 + 1)$.
 Here $\ln$ denotes Neperian logarithm.
We observe that the equality above is due to the fact that
$\msC(x)=\msC(-x)$ for all $x\in\IS^N$ and that
$\vol B_{\sin}(\bar{x}, \rho) = \vol B_{\IS}(\bar{x}, \theta)+
\vol B_{\IS}(-\bar{x}, \theta)$.

The same observation applies to the following result.

\begin{theorem}\label{thm:instance}
Let $\msC$ ba a conic condition number on $\R^{N+1}$ with set of ill-posed
inputs $\Sigma$. Assume that $\Sigma$ is contained in a real algebraic
hypersurface, given as the zero set of a homogeneous polynomial of
degree $d$. Let $\bar{x}\in \IS^N$ and $0\le \theta\le
\pi$. Then, for $\rho :=\sin\theta$,
$$
\Exp_{x\in B_{\IS}(\bar{x}, \theta)}\ln \msC(x)\le
\left\{\begin{array}{ll}
\ln \dfrac{Nd}{\rho + \frac{1-\rho}{\msC(\bar{x})}} +\ln 12 + 2&
  \mbox{if $\rho>\dfrac{1}{2\msC(\bar{x})+1}$}\\
\ln   \dfrac{1}{\rho + \frac{1-\rho}{\msC(\bar{x})}} +\ln 4&
  \mbox{if $\rho\leq\dfrac{1}{2\msC(\bar{x})+1}$.}
\end{array}\right.
$$
In particular, there is a uniform  explicit  bound $H(N,d,\theta, \msC(\bar x))$ --defined in \eqref{eq:UniformBound} below--  such that
$$  \Exp_{x\in B_{\IS}(\bar{x}, \theta)}\ln \msC(x)\le H(N,d,\theta,\msC(\bar x) ).$$
This bound satisfies  satisfies {\bf (LA0)}, since $H_\infty(N,d,\theta)= \ln \frac{Nd}{\sin\theta} + \Oh(1)$ as $H(N,d,\theta)$ in \eqref{eq:smooth-u},  {\bf (LA1)} and {\bf (LA2')}.
\end{theorem}

\proof
Assume first that
$\dfrac{1}{2\msC(\bar{x})+1}\le \rho\le 1$.
In this case, we have
$$
  \rho (2\msC(\bar{x})+1)\ge 1 \iff 2\rho \,
  \msC(\bar{x})+\rho \ge 1 \iff 2\rho\,\msC(\bar{x})
  \ge 1-\rho\iff  2\rho  \ge \frac{ 1-\rho }{ \msC(\bar{x})}
$$
and we can decompose
$$
  \rho = \frac{1}{3}\rho + \frac{1}{3}(2\rho) \ge
  \frac{1}{3}\Big(\rho +  \frac{1-\rho}{\msC(\bar{x})}\Big),
  \quad \mbox{i.e} \quad \frac{1}{\rho}\le
  \frac{3}{\rho + \frac{1-\rho}{\msC(\bar{x})}}.
$$
Therefore, by~\eqref{eq:smooth},
\begin{align*}
  \Exp_{x\in B_{\IS}(\bar{x}, \theta) } \ln \msC(x)&
  \le \ln \frac{Nd}{\rho} + \ln 4 + 2\\
  & \le \ln \frac{3Nd}{\rho + \frac{1-\rho}{\msC(\bar{x})} }
    + \ln 4 +2
  = \ln \frac{Nd}{\rho + \frac{1-\rho}{\msC(\bar{x})} }
  + \ln 12  +2.
\end{align*}
We next assume $0\le \rho < \dfrac{1}{2\msC(\bar{x})+1}$.
In this case,
$$
   \frac{1}{2\msC(\bar{x})+1} =
   \frac{1}{4}\Big(\frac{1}{2\msC(\bar{x})+1}+
   \frac{3}{2\msC(\bar{x})+1}\Big)
  > \frac{1}{4}\Big(\rho + \frac{1-\rho}{\msC(\bar{x})}\Big)
$$
since $\dfrac{3}{2 \msC(\bar{x})+1}>
\dfrac{1}{\msC(\bar{x})}>\dfrac{1-\rho}{\msC(\bar{x})}$.
Equivalently,
$$
   2 \msC(\bar{x})+1
  < \frac{4}{\rho +\frac{1-\rho}{\msC(\bar{x})}}.
$$
We also use here that for all
$x\in B_{\IS}(\bar{x},\theta)$,
\begin{align}\label{eq:small_rho}
\frac{1}{\msC(x)}& =\ d_{\sin}(x,\Sigma)\ \ge \
d_{\sin}(\bar{x},\Sigma)-d_{\sin}(x,\bar{x})\ \ge \
\frac{1}{\msC(\bar{x})} -\rho \nonumber \\ & \ \ge
\frac{1}{\msC(\bar{x})+\frac{1}{2}} -
\frac{1}{2\msC(\bar{x})+1}\ =\ \frac{1}{2\msC(\bar{x})+1},
\end{align}
and therefore
$$
   \msC(x)\le {2\msC(\bar{x})}+1
   < \frac{4}{\rho + \frac{1-\rho}{\msC(\bar{x})}}
%   \le \frac{4Nd}{\rho+ \frac{1-\rho}{\msC(\bar{x})}}
$$
which implies
\begin{align*}
  \ln  \msC(x) \le \ln \frac{1}{\rho + \frac{1-\rho}{\msC(\bar{x})}}
+\ln 4.
\end{align*}
This shows the first statement. We now derive the expression of a bound $H(N,d,\theta,\msC(\bar x))$.\\

Let $\varphi:[0,1]\to\R$ be the function
defined by
$$
     \rho\mapsto 2(Nd-1)
     \rho^{\log_{_{\!\!\frac{1}{2\msC(\bar{x})}}}\frac12} +1
$$
where  the exponent of $\rho$ in the numerator is
the logarithm in base $\dfrac{1}{2\msC(\bar{x}) }$ of
$\dfrac{1}{2}$, which, by continuity, we take to be 0 when
$\msC(\bar{x})=\infty$. We note that $\varphi$
is concave, monotonically increasing, d
satisfies $\varphi(0)=1$, $\varphi(1)=2Nd-1$, and when $\msC(\bar{x})=\infty$, $\varphi(\rho)=2Nd-1 $. Moreover, by monotonicity,
$$
   \varphi\Big(\frac{1}{2\msC(\bar{x})+1}\Big)
   \leq \varphi\Big(\frac{1}{2\msC(\bar{x})}\Big) =
   2(Nd-1)\frac12 +1= Nd.
$$
This implies, since $$\ln \frac{1}{\rho + \frac{1-\rho}{\msC(\bar{x})}}
+\ln 4 \le \ln \frac{\varphi(\rho)}
{\rho + \frac{1-\rho}{\msC(\bar{x})}} +\ln 12 +2\ \mbox{for} \ \ 0\le \rho < \dfrac{1}{2\msC(\bar{x})+1}$$ and using also  concavity,
$$\ln \frac{Nd}{\rho + \frac{1-\rho}{\msC(\bar{x})}}
+\ln 12+2\le \ln \frac{\varphi(\rho)}
{\rho + \frac{1-\rho}{\msC(\bar{x})}} +\ln 12+2 \ \mbox{for} \ \ \dfrac{1}{2\msC(\bar{x})+1}\le \rho \le 1,$$  that
$$
\Exp_{x\in B_{\IS}(\bar{x}, \theta) } \ln \msC(x)\le
\ln \frac{\varphi(\rho)}
{\rho + \frac{1-\rho}{\msC(\bar{x})}} +\ln 12 + 2.
$$
That is,
\begin{equation}\label{eq:UniformBound}H(N,d,\theta,\msC(\bar x))= \ln \frac{2(Nd-1)
     \rho^{\log_{_{\!\!\frac{1}{2\msC(\bar{x})}}}\frac12} +1}
{\rho + \frac{1-\rho}{\msC(\bar{x})}} +\ln 12 + 2. \end{equation}
Finally, it is trivial to verify, from the specific values taken by $\varphi$ mentioned previously,  that $ H(N,d,\theta,\msC(\bar x))$ satisfies {\bf (LA0)}, {\bf (LA1)} and {\bf (LA2')}.
\eproof

\section{The Gaussian case}\label{sec:gaussian}

We keep the same conic condition number $\msC$ but now
consider a Gaussian measure $N(\bar{x},\sigma^2\Id)$
in $\R^{N+1}$ centered at
$\bar{x}\in\IS^N$ and with covariance matrix
$\sigma^2\Id$ for $0<\sigma <\infty$, that is with density function
given by
$$
    \frac{1}{(2\pi\sigma^2)^{\frac{N+1}{2}}}\exp
    \Big(\frac{-\|x-\bar x\|^2}{2\sigma^2}\Big).
$$
Since our local analysis will rely on a smoothed analysis in this case, which is not yet known, we begin by studying a general smoothed analysis for the Gaussian case.

\subsection{Smoothed analysis}

Let $\bar x \in \IS^N$. We recall that, for any $0\le \theta\le \frac{\pi}{2}$,
$$
 B_{\IS}(\bar{x},\theta)=
 \{  x\in \IS^N\,:\, 0\le\sphericalangle(x,\bar x) < \theta \} ,
$$
and in the particular case  $\theta =\frac{\pi}{2}$ we denote
$$
  \IS^N_+({\bar x}):=B_{\IS}\Big(\bar x,\frac{\pi}{2}\Big)
  = \Big\{x\in \IS^N\,:\, 0\le \sphericalangle(x,\bar x) < \frac{\pi}{2}\Big\}
  =\big\{x\in \IS^N\,:\, \langle x, \bar x\rangle > 0\big\},
$$
the open half-sphere centered at $\bar x$.\\

The main result of this section is the following smoothed analysis for the Gaussian distribution.
\begin{theorem}\label{th:smoothedGaussian}
Let $\msC$ be a conic condition number on $\R^{N+1}$ with set of ill-posed
inputs $\Sigma$. Assume that $\Sigma$ is contained in a real algebraic
hypersurface, given as the zero set of a homogeneous polynomial of
degree $d$, and that $N\ge 5$.
 Then, there exists an explicit bound $H(N,d,\sigma)$ --defined in \eqref{SmoothGaussianBound}--
 such that
$$\max_{\bar x\in \IS^N}
 \Exp_{x\sim N(\bar{x}, \sigma^2\Id)}\ln \msC(x)
 \le H(N,d,\sigma).
 $$
 This bound satisfies
 \begin{description}
\item[(SA1)]  $\displaystyle \lim_{\sigma\to 0}H(N,d,\sigma)=\infty$, the worst-case value.
\item[(SA2)] $\displaystyle \lim_{\sigma\to \infty}H(N,d,\sigma)=\ln(Nd)+2(\ln 2+1)$,    the average value,
in remarkable coincidence with~\eqref{eq:smooth-pi}.
 \end{description}

\end{theorem}

The following map
plays a central role in all what follows,
\begin{equation}\label{eq:Psi}
  \Psi: \R^{N+1}\setminus \bar{x}^\perp \to \IS^N_+(\bar{x}),
  \qquad
  x\mapsto \begin{cases} \|x\|^{-1}x& \mbox{if } \langle x,\bar x\rangle > 0\\
 -\|x\|^{-1}x& \mbox{otherwise.} \end{cases}
\end{equation}

The main stepping stone towards the proof of
Theorem~\ref{th:smoothedGaussian} is the following.

\begin{proposition}\label{prop:tech}
Let $\bar x \in \IS^N$. There exists a probability density
$f: [0,\frac{\pi}{2}]\to \R_{\ge 0}$ of a random variable
$\theta \in [0,\frac{\pi}{2}]$, associated to $\bar x$, $\sigma$ and $N$, such that
for all measurable function $F:\R^{N+1}\to \R_{\ge 0}$  satisfying
$F(x)=F(\lambda x)$ for all $\lambda\in \R^\times$, one has
$$
  \Exp_{y\sim N(\bar{x},\sigma^2\Id)} F(y)
  = (1-e^{-\frac{1}{2\sigma^2}})\Exp_{\theta \sim f} \Big(\Exp_{x\in B_{\IS}(\bar{x},\theta)} F(x)\Big)  \  + \  e^{-\frac{1}{2\sigma^2}}\Exp_{x\in \IS^N_+(\bar x)} F(x).
$$
\end{proposition}

We begin by proving the following lemma.

\begin{lemma}\label{lem:tech}
For any measurable function $F:\R^{N+1}\to \R_+$ satisfying
$F(\lambda y)=F(y), \ \forall \lambda \in \R^\times$, one has
$$
  \Exp_{y\sim N(\bar{x}, \sigma^2\Id)} F(y)=
  \int_{\IS^N_+(\bar x)} G_{\bar x}(\sphericalangle(x,\bar x)) F(x)
  \diff x
$$
where $G_{\bar x}: [0,\frac{\pi}{2}] \to \R_{>0}$ is a decreasing function of
$\alpha$ defined by
\begin{align*}
   G_{\bar x}(\alpha)&=\frac{1}{(2\pi\sigma^2)^{\frac{N+1}{2}}}
   \int_{-\infty}^{\infty}
  \exp\Big(-\frac{\lambda^2+1 - 2 \lambda \cos\alpha  }{2\sigma^2}\Big)
  |\lambda |^N\diff \lambda.
\end{align*}
\end{lemma}

\proof
We have
\begin{align*}
\Exp_{y\sim N(\bar{x}, \sigma^2\Id)} F(y)
&= \frac{1}{(2\pi\sigma^2)^{\frac{N+1}{2}}}
\int_{\R^{N+1}} F(y)
\exp\Big(\frac{-\|y-\bar{x}\|^2}{2\sigma^2}\Big)
\diff y\nonumber\\
&= \frac{1}{(2\pi\sigma^2)^{\frac{N+1}{2}}}\int_{\IS^N_+(\bar x)} \Big(\int_{-\infty}^{\infty}
F(\lambda x)
\exp\Big(\frac{-\|\lambda x-\bar{x}\|^2 }{2\sigma^2}\Big)|\lambda|^N
\diff \lambda\Big)\diff x \nonumber\\
&=\int_{\IS^N_+(\bar x)} F(x) \bigg[\frac{1}{(2\pi\sigma^2)^{\frac{N+1}{2}}}\int_{-\infty}^{\infty}
\exp\Big(\frac{-\|\lambda x-\bar{x}\|^2 }{2\sigma^2}\Big)
|\lambda|^N\diff \lambda\bigg]\diff x\\
& = \int_{\IS^N_+(\bar x)} F(x) G(x)\diff x
\end{align*}
where the second equality follows from  the transformation
formula~\cite[Thm.~2.1]{Condition} applied to the diffeomorphism
$$
 \Phi: \R^{N+1}\setminus \bar{x}^\perp \to \IS^N_+(\bar{x})\times \R\setminus\{0\},
 \qquad
 x\mapsto \begin{cases}
   (\Psi(x),\|x\|^2)& \mbox{if }
  \langle x,\bar x\rangle > 0\\
  (\Psi(x),-\|x\|^2)& \mbox{otherwise,} \end{cases}
$$
and
$$
  G(x):= \frac{1}{(2\pi\sigma^2)^{\frac{N+1}{2}}}\int_{-\infty}^{\infty}
  \exp\Big(\frac{-\|\lambda x-\bar{x}\|^2 }{2\sigma^2}\Big)
  |\lambda|^N\diff \lambda
$$
does not depend on $F$.
Now, for $\bar x, x\in \IS^N_+(\bar x)$,
$$\|\lambda x-\bar{x}\|^2=\lambda^2-2\lambda \cos(\sphericalangle(x,\bar x)) +1.
$$
Therefore, $G(x)=:G_{\bar x}(\sphericalangle(x,\bar x))$ where
for $0\le \alpha \le \frac{\pi}{2}$,
$$
  G_{\bar x}(\alpha)= \frac{1}{(2\pi\sigma^2)^{\frac{N+1}{2}}}\int_{-\infty}^{\infty}
  \exp\Big( - \frac{\lambda^2+1 - 2 \lambda \cos\alpha }{2\sigma^2}\Big)
  |\lambda|^N\diff \lambda,
$$
which is a continuously differentiable decreasing function of $\alpha$.
\eproof

\proofof{Proposition~\ref{prop:tech}}
By Lemma~\ref{lem:tech},
\begin{equation}\label{eq1}
  \Exp_{y\sim N(\bar{x}, \sigma^2\Id)} F(y)=
 \int_{\IS^N_+(\bar x)} G_{\bar x}(\sphericalangle(x,\bar x) ) F(x)
  \diff x.
\end{equation}
Now, by the fundamental Theorem of Calculus for $0<\alpha<\frac{\pi}{2}$,
$$
G_{\bar x}(\alpha)= G_{\bar x}\Big(\frac{\pi}{2}\Big)
-\int_\alpha^{\frac{\pi}{2}}G'_{\bar x}(\theta)\diff \theta=
G_{\bar x}\Big(\frac{\pi}{2}\Big)
-\int_0^{\frac{\pi}{2}}\ \uno_{\{\alpha\le \theta\}}G'_{\bar x}(\theta)\diff \theta.
$$
Replacing this in \eqref{eq1} and changing the order of integration, we
obtain
$$
  \Exp_{y\sim N(\bar{x}, \sigma^2\Id)} F(y)= G_{\bar x}\Big(\frac{\pi}{2}\Big)
  \int_{\IS^N_+(\bar x)}  F(x)
  \diff x - \int_0^{\frac{\pi}{2}} \Big(\int_{\IS^N_+(\bar x)}F(x)
  \uno_{\{\sphericalangle(x,\bar x)\le \theta\}}dx\Big)G'_{\bar x}(\theta)\diff \theta.
$$
Now, since
$$
   \Exp_{x\in \IS^N_+(\bar x)}F(x)=
   \frac{\displaystyle \int_{\IS^N_+(\bar x)}F(x)\diff x}{\vol(\IS^N_+(\bar x))}
\quad \mbox{
and
} \quad
  \Exp_{x\in B_\IS(\bar x,\theta)}F(x)
 = \frac{\displaystyle \int_{B_\IS(\bar x,\theta)}F(x)\diff x}
  {\vol(B_{\IS}(\bar{x},\theta))},
$$
we obtain
\begin{align*}
  \Exp_{y\sim N(\bar{x},\sigma^2\Id)} F(y)
  =\;&G_{\bar x}\Big(\frac{\pi}{2}\Big)
  \vol(\IS^N_+(\bar x))\Exp_{x\in \IS^N_+(\bar x)}F(x)\\
  &- \int_0^{\frac{\pi}{2}} \Big(\vol(B_{\IS}(\bar{x},\theta))
  \Exp_{x\in B_\IS(\bar x,\theta)}F(x)\Big)G'_{\bar x}(\theta)\diff \theta.
\end{align*}
We now denote $$f(\theta):=- \frac{\vol(B_{\IS}(\bar{x},\theta))G'_{\bar x}(\theta)}{1-e^{\frac{-1}{2\sigma^2}}},$$
which is a non-negative function since $G_{\bar x}$ is decreasing,
and rewrite the equality above as
\begin{equation*}
 \Exp_{y\sim N(\bar{x},\sigma^2\Id)} F(y)=H(N,\sigma)\Exp_{x\in \IS^N_+(\bar x)}F(x)
 +(1-e^{\frac{-1}{2\sigma^2}})\int_0^{\frac{\pi}{2}}\Big(\Exp_{x\in B_{\IS}(\bar{x},\theta)} F(x)\Big)
 f(\theta)\diff \theta ,
\end{equation*}
where
$$
  H(N,\sigma)=G_{\bar x}\Big(\frac{\pi}{2}\Big) \vol(\IS^N_+(\bar x))
  = \vol(\IS^N_+(\bar x))\frac{1}{(2\pi\sigma^2)^{\frac{N+1}{2}}}\int_{-\infty}^{\infty}
  \exp\Big( - \frac{\lambda^2+1 }{2\sigma^2}\Big)
  |\lambda|^N\diff \lambda.
$$
We now prove that $H(N,\sigma)= e^{-\frac{1}{2\sigma^2}}$:
Changing variables $\nu=\frac{\lambda}{\sigma}$ we have
\begin{eqnarray*}
  H(N,\sigma)&=&\frac{ \vol(\IS^N_+(\bar x))}{(2\pi\sigma^2)^{\frac{N+1}{2}}}
  \int_{-\infty}^{\infty}
  \exp\Big( - \frac{\lambda^2+1 }{2\sigma^2}\Big)
  |\lambda|^N\diff \lambda\\
  &=& \frac{ \vol(\IS^N_+(\bar x))}{(2\pi\sigma^2)^{\frac{N+1}{2}}}
  \int_{-\infty}^{\infty}
  \exp\Big( - \frac{\nu^2}{2}-\frac{1 }{2\sigma^2}\Big)
  |\nu|^N \sigma^{N+1} \diff \nu\\
  &=& e^{-\frac{1 }{2\sigma^2}}
  \bigg[\frac{\vol(\IS^N_+(\bar x))}{(2\pi)^{\frac{N+1}{2}}}
  \int_{-\infty}^{\infty}
  e^{- \frac{\nu^2}{2}}
  |\nu|^N \diff \nu\bigg].
\end{eqnarray*}
To estimate the quantity between the square brackets we use
the known equality
$$
 \int_0^\infty \nu^{N}e^{-\frac{x^2}{2}}\diff \nu
 = \Gamma\Big(\frac{N+1}{2}\Big)2^{\frac{N-1}{2}}
$$
together with~\eqref{eq:volSN} to obtain
\begin{eqnarray*}
   \frac{\vol(\IS^N_+(\bar x))}{(2\pi)^{\frac{N+1}{2}}}
  \int_{-\infty}^{\infty}
  \exp\Big(- \frac{\nu^2}{2}\Big)
  |\nu|^N \diff \nu &=&\frac{\vol(\IS^N_+(\bar x))}{(2\pi)^{\frac{N+1}{2}}}
  \Gamma\Big(\frac{N+1}{2}\Big)2^{\frac{N+1}{2}}\\
  &=& \frac{\pi^{\frac{N+1}{2}}}{\Gamma\Big(\frac{N+1}{2}\Big)(2\pi)^{\frac{N+1}{2}}}
  \Gamma\Big(\frac{N+1}{2}\Big)2^{\frac{N+1}{2}}\\
  &=& 1.
\end{eqnarray*}
Therefore
\begin{equation*}
 \Exp_{y\sim N(\bar{x},\sigma^2\Id)} F(y)=e^{\frac{-1}{2\sigma^2}}\Exp_{x\in \IS^N_+(\bar x)}F(x)
 +(1-e^{\frac{-1}{2\sigma^2}})\int_0^{\frac{\pi}{2}}\Big(\Exp_{x\in B_{\IS}(\bar{x},\theta)} F(x)\Big)
 f(\theta)\diff \theta.
\end{equation*}
This implies, by taking $F=1$, that
$$1 = e^{\frac{-1}{2\sigma^2}} + (1-e^{\frac{-1}{2\sigma^2}})\int_0^{\frac{\pi}{2}}
 f(\theta)\diff \theta,$$
 i.e.
 $$\int_0^{\frac{\pi}{2}}
 f(\theta)\diff \theta=1.$$
Therefore  $f$ is a density on $[0,\frac{\pi}{2}]$, and
$$\int_0^{\frac{\pi}{2}}\Big(\Exp_{x\in B_{\IS}(\bar{x},\theta)} F(x)\Big)
 f(\theta)\diff \theta = \Exp_{\theta\sim f} \Big( \Exp_{x\in B_{\IS}(\bar{x},\theta)} F(x) \Big) .$$
\eproof

\bigskip
Since $\msC(x)=\msC(\lambda x)$ for all $\lambda \in \R^\times$, we can
now focus on $F(x):=\ln\msC(x)$.

\begin{proposition} \label{together}
With the notation in Proposition~\ref{prop:tech}, we have
\begin{equation*}
\Exp_{y\sim N(\bar{x},\sigma^2\Id)} \ln \msC(y)
\le   (1-e^{\frac{-1}{2\sigma^2}})\Exp_{\theta \sim f}(\ln \Big( \frac{1}{\sin\theta}\Big)) + \ln (Nd)+ 2(\ln 2 + 1).
\end{equation*}
\end{proposition}

\proof
Replacing the expectations in the right-hand side of
the equality in Proposition~\ref{prop:tech} by their bound
in~\eqref{eq:smooth} for   $\rho=\sin \theta$ and $\rho=1=\sin\frac{\pi}{2}$,
we obtain
\begin{eqnarray*}
 \Exp_{y\sim N(\bar{x},\sigma^2\Id)} \ln\msC(y)
 &=& (1-e^{\frac{-1}{2\sigma^2}}) \Exp_{\theta \sim f} \Big( \Exp_{x\in B_{\IS}(\bar{x},\theta)} \ln\msC(x) \Big)+e^{-\frac{1}{2\sigma^2}}
 \Exp_{x\in {\IS_+^N(\bar x)}}\ln \msC(x)\\
 &\le &(1-e^{\frac{-1}{2\sigma^2}}) \Exp_{\theta \sim f}(\ln \frac{Nd}{\sin\theta} + K)
      +e^{-\frac{1}{2\sigma^2}}
 \Exp_{x\in {\IS_+^N(\bar x)}}(\ln (Nd) + K)\\
 & \le & (1-e^{\frac{-1}{2\sigma^2}}) \Big(\Exp_{\theta \sim f}(\ln \Big( \frac{1}{\sin\theta}\Big) +
   \ln (Nd) + K) \Big)   +e^{-\frac{1}{2\sigma^2}}
(\ln (Nd) + K)\\
&\le & (1-e^{\frac{-1}{2\sigma^2}})\Exp_{\theta \sim f}(\ln \Big( \frac{1}{\sin\theta}\Big)) + \ln (Nd) + K,
\end{eqnarray*}
where $K=2(\ln 2 + 1)$. The result follows from the last equality
in Proposition~\ref{prop:tech}.
\eproof

Our next goal is to estimate the right-hand side in
Proposition~\ref{together}.

\begin{lemma}\label{lem:A2}
Let $0\le t \le \frac{\pi}{4}$. Then
\begin{align*}
\int_{t}^{\frac{\pi}{2}}\ln\Big(\frac1{\sin\theta}\Big)f(\theta)
 \diff \theta
 \le  \ln  \sqrt 2 + \int_{\sin t}^{\frac{\sqrt 2}{2}}
 \Big(\int_{t}^{\arcsin s}
 f(\theta)\diff \theta\Big)\frac1s\diff s.
\end{align*}
\end{lemma}

\proof
Write
$$
\int_{t}^{\frac{\pi}{2}}\ln\Big(\frac1{\sin\theta}\Big)f(\theta)
 \diff \theta =
\int_{t}^{\frac{\pi}{4}}\ln\Big(\frac1{\sin\theta}\Big)f(\theta)
 \diff \theta+
\int_{\frac{\pi}{4}}^{\frac{\pi}{2}}\ln\Big(\frac1{\sin\theta}\Big)f(\theta)
 \diff \theta.$$
Since $\frac{1}{\sin\theta}\le \sqrt 2$ for $\theta \in
[\frac{\pi}{4},\frac{\pi}{2}]$, the second term satisfies
$$\int_{\frac{\pi}{4}}^{\frac{\pi}{2}}\ln\Big(\frac1{\sin\theta}\Big)f(\theta)
\diff \theta \le \ln \sqrt 2 \int_{\frac{\pi}{4}}^{\frac{\pi}{2}}f(\theta)
\diff \theta.$$
We analyze the first term.
Let $A=\{ (\theta,r) \in [t,\frac{\pi}{4}]\times
[0,\ln\big(\frac{1}{\sin t}\big)]\,:\,
r\le \ln (\frac1{\sin\theta})\}$, and
$A_r=\{\theta \in [t,\frac{\pi}{4}]: r\le \ln (\frac1{\sin\theta})\}$.
By Fubini's Theorem we have both
$$
 \int_{(\theta, r)\in A} f(\theta)\diff(\theta,r)
 = \int_t^{\frac{\pi}{4}}\Big( \int_{0}^{ \ln (\frac1{\sin\theta} )}
 \diff r\Big) f(\theta)\diff \theta
 = \int_t^{\frac{\pi}{4}}\ln\Big(\frac1{\sin\theta}\Big)f(\theta)\diff\theta
 $$
and
\begin{align*}
 \int_{(\theta, r)\in A} f(\theta)\diff(\theta,t)
  & \ = \ \int_{0}^{\ln (\frac{1}{\sin t})} \Big(\int_{\theta\in A_r} f(\theta)
  \diff \theta\Big) \diff r\\
   & \ = \ \int_{0}^{\ln\sqrt 2} \Big(\int_{\theta\in A_r} f(\theta)
  \diff \theta\Big) \diff r +
  \int_{\ln\sqrt2}^{\ln (\frac{1}{\sin t})} \Big(\int_{\theta\in A_r} f(\theta)
  \diff \theta\Big) \diff r\\
 & = \
  \ln \sqrt 2 \int_t^{\frac{\pi}{4}} f(\theta)\diff \theta +
  \int_{\ln \sqrt 2}^{\ln (\frac{1}{\sin t}) } \Big(\int_{\sin t\le \sin \theta \le e^{- r}}
  f(\theta)\diff \theta\Big) \diff r,
\end{align*}
since $t\le \frac{\pi}{4}$ implies $\ln \sqrt 2\le \ln\big(\frac{1}{\sin t}\big) $
and when $r\le \ln \sqrt 2$, then $A_r=[t,\frac{\pi}{4}]$. Therefore,
$$
\int_{(\theta, r)\in A} f(\theta)\diff(\theta,r)
\ = \ \ln \sqrt 2 \int_t^{\frac{\pi}{4}} f(\theta)\diff \theta
+ \int_{\sin t}^{\frac{\sqrt 2}{2}} \Big(\int_{t}^{\arcsin s} f(\theta)
  \diff \theta\Big)\frac1s\diff s ,
$$
by taking $s=e^{-r}$.
Finally,
\begin{align*}
\int_{t}^{\frac{\pi}{2}}\ln\Big(\frac1{\sin\theta}\Big)f(\theta)
\diff \theta & \le  \ln \sqrt 2 \int_{\frac{\pi}{4}}^{\frac{\pi}{2}}f(\theta)
+   \ln \sqrt 2 \int_t^{\frac{\pi}{4}} f(\theta)\diff \theta
+ \int_{\sin t}^{\frac{\sqrt 2}{2}} \Big(\int_{t}^{\arcsin s} f(\theta)
\diff \theta\Big)\frac1s\diff s \\
& \le \ln \sqrt 2 \int_t^{\frac{\pi}{2}} f(\theta)\diff \theta
+ \int_{\sin t}^{\frac{\sqrt 2}{2}} \Big(\int_{t}^{\arcsin s} f(\theta)
\diff \theta\Big)\frac1s\diff s\\
& \le  \ln \sqrt 2
  + \int_{\sin t}^{\frac{\sqrt 2}{2}} \Big(\int_{t}^{\arcsin s} f(\theta)
  \diff \theta\Big)\frac1s\diff s
\end{align*}
  since $\int_t^{\frac{\pi}{2}}f(\theta)\diff \theta \le 1$.
\eproof

\begin{lemma}\label{lem:B}
Assume $N\ge 5$.
For all $t\in[0,\frac{\pi}{4}]$, one has
$$
  \int_{0}^ t f(\theta)\diff \theta \leq
  \min\Big\{1,\frac{1}{1-e^{\frac{-1}{2\sigma^2}}}\Big(\frac12 (\sin(2t))^N +
  \frac{(\sin t)^N}{\sigma^{N+1}}\Big)\Big\}.
$$
\end{lemma}

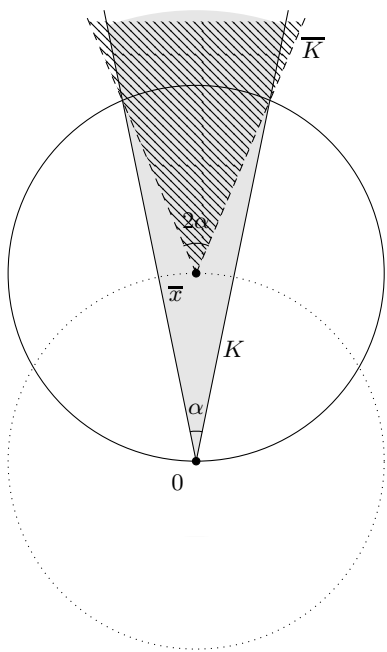
\begin{figure}[H]\centering
\begin{tikzpicture}[scale=.5,point/.style={draw,minimum size=0pt,
	inner sep=1pt,circle,fill=black}]
\begin{scope}
\fill[gray!20!white] (0,0) circle (7);
\fill[white] (-0.4,-7)-- (2.45,7) -- (7,7) -- (7,-7) -- (-0.4,-7);
\fill[white] (0.4,-7)-- (-2.45,7) -- (-7,7) -- (-7,-7) -- (0.4,-7);
\fill[pattern=north west lines] (0,0)-- (-2.85,6.7) -- (2.85,6.7) -- (0,0);
\end{scope}
\draw (0,0) node(x) [point,label=220:{\footnotesize$\bar{x}$}] {};
\draw(x) circle [radius=5];
\draw[dotted](0,-5) circle [radius=5];
\draw (0,-5) node(q) [point,label=-120:{\footnotesize$0$}] {};
\draw[dashed] (0,0) --(2.9,6.8);
\draw[dashed] (0,0) --(-2.9,6.8);
\draw (q) -- (2.45,7);
\draw (q) -- (-2.45,7);
\draw ([shift={(0,0)}]65:0.8) arc [radius =0.8,
start angle=65, end angle=115];
\path (0,1.4) node{\footnotesize $2\alpha$};
\draw ([shift={(0,-5)}]78:0.8) arc [radius =0.8,
start angle=78, end angle=102];
\path (0,-3.6) node{\footnotesize $\alpha$};
\path (1,-2) node{\footnotesize $K$};
\path (3.1,6) node{\footnotesize $\bar{K}$};
\end{tikzpicture}
\caption{{\small The cones $K$ (shaded) and $\bar{K}$
(line patterned).}}
\label{fig:cones}
\end{figure}

\proof
For $t\le \frac{\pi}{2}$,
\begin{align*}
  \int_{0} ^{ t} f(\theta)\diff \theta
  &= \Exp_{\theta\sim f}(\uno_{\{\theta\leq t\}})
  \;\le\; \Exp_{\theta\sim f} \big(\Exp_{B_{\IS}(\bar{x},\theta)}
  \uno_{\{\sphericalangle(x,\bar{x})\leq t\}}
  \big) )\\
  &\leq \frac{1}{1-e^{\frac{-1}{2\sigma^2}}}\Exp_{y\sim N(\bar{x}, \sigma^2\Id)}
  (\uno_{ \{\sphericalangle (\Psi(y),\bar{x} )\leq t \}})\\
&=\frac{1}{1-e^{\frac{-1}{2\sigma^2}}} \Prob_{y\sim N(\bar{x}, \sigma^2\Id)}
     \big\{\sphericalangle (\Psi (y),\bar{x} )\leq t\big\},
\end{align*}
for $\Psi$ defined in \eqref{eq:Psi}.
The first inequality holds because for $\theta \le t$,
$\sphericalangle(x,\bar{x})\leq \theta $ implies
$\sphericalangle(x,\bar{x})\leq t $, and the second by
Proposition~\ref{prop:tech} applied to
$F=\uno_{\big\{\sphericalangle (\Psi(y),\bar{x})\le t\big\}}$.  It is
then enough to bound the right-hand expression.\\ We observe that for
$0\le t\le \frac{\pi}{2}$, the set
$K=\big\{y\in\R^{N+1}\,:\,\sphericalangle(\Psi(y),\bar{x})\leq
t\big\}$ is a pointed cone with vertex at $0$, central axis passing
through $\bar{x}$ and angular opening $\alpha:=2t$.
In addition, one can prove by the cosine theorem that this cone is
included in the union of the pointed cone $\bar{K}$ with vertex at
$\bar{x}$, central axis passing through $2\bar{x}$ and angular opening
$2\alpha$ with the intersection $K\cap B(\bar{x},1)$ (see
Figure~\ref{fig:cones}). Hence, the measure of $K$ (with respect to
$N(\bar{x},\sigma^2\Id)$) is bounded by the sum of the measures of
$\bar{K}$ and $K\cap B(\bar{x},1)$.
As the vertex $\bar{x}$ of $\bar{K}$ coincides with the center of
$N(\bar{x},\sigma)$, the measure of $\bar{K}$ with respect to $N(\bar{x},\sigma)$
equals the proportion of the volume (in $\IS(\bar{x},1)$) of the intersection of
$\bar{K}$ with $\IS(\bar{x},1)$ within this sphere. That is,
the measure of $\bar{K}$ with respect to $N(\bar{x},\sigma)$ satisfies
$$ \Prob_{x\sim N(\bar{x},\sigma^2\Id)}\{x\in \bar{K}\} =
\frac{\vol(B_{\IS}(\bar{x},2t))}{\Oh_{N}} $$ where, we recall,
$\Oh_N:=\vol(\IS^N)$.
Using~\eqref{eq:volBS} we deduce that, for $t\in[0,\frac{\pi}{4}]$,
$$
   \Prob_{x\sim N(\bar{x},\sigma^2\Id)}\{x\in \bar{K}\} \leq \frac12
   (\sin(2t))^N.
$$
Also,
\begin{align*}
 \Prob_{x\sim N(\bar{x},\sigma^2\Id)}&\{x\in K\cap B(\bar{x},1)\}
\ = \ \int_{x\in K\cap B(\bar{x},1)}
  \frac{1}{(2\pi\sigma^2)^{\frac{N+1}{2}}}
 \exp\Big(-\frac{\|x-\bar{x}\|^2}{2\sigma^2}\Big)\diff x \\
 &\leq \frac{1}{(2\pi\sigma^2)^{\frac{N+1}{2}}}
 \int_{x\in K\cap B(\bar{x},1)} 1\diff x \
= \ \frac{1}{(2\pi\sigma^2)^{\frac{N+1}{2}}}
 \vol(K\cap B(\bar{x},1))\\
 &\le\frac{1}{(2\pi\sigma^2)^{\frac{N+1}{2}}}\vol(K\cap B(0,2))
\ = \ \frac{\vol(B_\IS(\bar x, t))}{(2\pi\sigma^2)^{\frac{N+1}{2}}\Oh_N}
 \vol(B(0,2))\\
&\underset{\eqref{eq:volBS}}{\le}
  \frac{2^{N+1}(\sin t)^N}{(2\pi\sigma^2)^{\frac{N+1}{2}}\cdot 2}
 \vol(B(0,1))\ \underset{\eqref{eq:volSN}\eqref{eq:volB}}{\le} \
 \frac{2^{\frac{N+1}{2}}(\sin t)^N}{\Gamma(\frac{N+1}{2})(N+1)\sigma^{N+1}}
\\	
&\le \frac{2^{N+\frac12}e^{\frac{N-1}{2}}(\sin t)^N}{\sqrt{\pi}
(N-1)^{\frac{N}{2}}(N+1)\sigma^{N+1}}.
\end{align*}
Here we used the well-known lower bound
$\Gamma(\frac{N+1}{2})>\sqrt{2\pi}\Big(\frac{N-1}{2}\Big)^{\frac{N}{2}}
e^{-\frac{N-1}{2}}$ (see for instance~\cite[Eq.~2.14]{Condition})
for the last inequality. We finish the proof by noting that it can be
easily proven by induction, using for instance that $N^{N+1}\ge 2N(N-1)^N$,
that for all
$N\ge 5$,  we have
\begin{equation}\tag*{\qed}
  \frac{2^{N+\frac12}e^{\frac{N-1}{2}}}{\sqrt{\pi}
  (N-1)^{\frac{N}{2}}(N+1)}\leq 1.
\end{equation}

\begin{lemma}\label{lem:int1}
Assume   $N\ge 5$. Then,
$$
\Exp_{\theta\sim f}(\ln\Big(\frac1{\sin\theta}\Big))\le
 \frac{1}{N} \Big(1+ \ln \big( 2^{N-1}+\frac{1}{\sigma^{N+1}}\big)
 -  \ln
(1-e^{-\frac{1}{2\sigma^2}}) \Big).
$$
\end{lemma}

\proof
We have by Lemma~\ref{lem:A2} with $t=0$,
$$
\Exp_{\theta\sim f}(\ln\Big(\frac1{\sin\theta}\Big))
 \le    \ln  \sqrt 2 + \int_0^{\frac{\sqrt 2}{2}}
 \Big(\int_{0}^{\arcsin s}
 f(\theta)\diff \theta\Big)\frac1s\diff s,
$$
where by Lemma~\ref{lem:B}, since $0\le \arcsin s \le \frac{\pi}{4}$ for
$0\le s\le \frac{\sqrt 2}{2}$,
\begin{align*}
  \int_{0}^{\arcsin s} f(\theta)\diff \theta & \leq
  \min\bigg\{1, \frac{1}{1-e^{-\frac{1}{2\sigma^2}}}\Big(\frac12 (\sin(2\arcsin s))^N +
  \frac{(\sin(\arcsin s))^N}{\sigma^{N+1}}\Big)\bigg\}
  \\
  & \le \min\bigg\{1, \frac{1}{1-e^{-\frac{1}{2\sigma^2}}}\Big(2^{N-1}s^N +
  \frac{s^N}{\sigma^{N+1}}\Big)\bigg\}
  \\
  & \le \min\Big\{1, \frac{2^{N-1} +
  \frac{1}{\sigma^{N+1}}}{1-e^{-\frac{1}{2\sigma^2}}}s^N\Big\}.
\end{align*}
We have
$$\frac{2^{N-1} +
  \frac{1}{\sigma^{N+1}}}{1-e^{-\frac{1}{2\sigma^2}}}s^N\le 1 \iff \Big(2^{N-1} +
  \frac{1}{\sigma^{N+1}}\Big)s^N \le 1-e^{-\frac{1}{2\sigma^2}}  \iff
  s\le c(N,\sigma),
$$
where  $c(N,\sigma):=\sqrt[N]
{\dfrac{(1-e^{-\frac{1}{2\sigma^2}})\sigma^{N+1}}{1+2^{N-1}\sigma^{N+1}}}$.
In addition we observe that for all
$N\ge 2$, $c(N,\sigma)< \dfrac{\sqrt 2}{2}$ since
\begin{eqnarray*}
c(N,\sigma) <\frac{1}{\sqrt{2}}
&\iff & \dfrac{(1-e^{-\frac{1}{2\sigma^2}})\sigma^{N+1}}{1+2^{N-1}\sigma^{N+1}}
          <\frac{1}{2^\frac{N}{2}}\\
&\iff & 2^\frac{N}{2}(1-e^{-\frac{1}{2\sigma^2}}) \sigma^{N+1} < 1+2^{N-1}\sigma^{N+1}.
\end{eqnarray*}
Rewriting $c(N,\sigma)^{-N}= \dfrac{2^{N-1}+\frac{1}{ \sigma^{N+1}}}{ 1-e^{-\frac{1}{2\sigma^2}}}
$ we get
\begin{align*}
\Exp_{\theta\sim f}(\ln\Big(\frac1{\sin\theta}\Big))
 &\le   \ln \sqrt 2
 +  \int_0^{c(N,\sigma)}\Big(\int_{0}^{\arcsin s}
 f(\theta)\diff \theta\Big)\frac1s\diff s
+ \int_{c(N,\sigma)}^{\frac{\sqrt 2}{2}}\Big(\int_{0}^{\arcsin s}
 f(\theta)\diff \theta\Big)\frac1s\diff s\\
 &\le   \ln \sqrt 2 +
 \int_0^{c(N,\sigma)}c(N,\sigma)^{-N}s^{N-1}\diff s
  + \int_{c(N,\sigma)}^{\frac{\sqrt 2}{2}}\frac{1}{s}\diff s\\
  &\le  \ln \sqrt 2+ \frac{1}{N}
  + \ln \frac{\sqrt 2}{2}- \ln c(N,\sigma) \\
& = \frac{1}{N} \Big( 1 +  \ln \frac{2^{N-1}+\frac{1}{\sigma^{N+1}}}
{1-e^{-\frac{1}{2\sigma^2}}} \Big)\\
& = \frac{1}{N}
 \Big(1+ \ln \big( 2^{N-1}+\frac{1}{\sigma^{N+1}}\big)
 - \ln(1-e^{-\frac{1}{2\sigma^2}}) \Big).
\end{align*}
\eproof

\proofof{Theorem~\ref{th:smoothedGaussian}}
By Proposition~\ref{together} and Lemma~\ref{lem:int1},
\begin{align*}
\Exp_{y\sim N(\bar{x},\sigma^2\Id)} \ln \msC(y)
&\le   (1-e^{\frac{-1}{2\sigma^2}})\Exp_{\theta \sim f}(\ln \Big( \frac{1}{\sin\theta}\Big)) + \ln (Nd) +K \\
 & \le
\frac{1-e^{\frac{-1}{2\sigma^2}}}{N} \Big(1+ \ln \big( 2^{N-1}+\frac{1}{\sigma^{N+1}}\big)
 - \ln(1-e^{-\frac{1}{2\sigma^2}}) \Big)\\
 & \quad +\ln(Nd) + K,
\end{align*}
with $K=2(\ln 2 + 1)$.
We then define
\begin{equation}\label{SmoothGaussianBound}
H(N,d,\sigma)= \frac{(1-e^{-\frac{1}{2\sigma^2}})}{N}
 \Big(1+ \ln \big( 2^{N-1}+\frac{1}{\sigma^{N+1}}\big) -  \ln
 (1-e^{-\frac{1}{2\sigma^2}}) \Big)+ \ln (Nd)+2(\ln 2 +1).
\end{equation}
We now verify that $H(N,d,\sigma)$ satisfies {\bf (SA1)} and {\bf (SA2)}:
\begin{description}
\item[(SA1)] $\displaystyle \lim_{\sigma\to 0} H(N,d,\sigma)= \displaystyle\lim_{\sigma\to 0}\Big(\dfrac{1}{N}\Big(1+ \ln \big( 2^{N-1}+\dfrac{1}{\sigma^{N+1}}\big) \Big)+ \ln (Nd)
+ 2(\ln 2 +1)\Big)$\\
${\ } \qquad \qquad \qquad  \quad =\displaystyle \lim_{\sigma\to \infty} \Big(\dfrac{N+1}{N}\ln\dfrac{Nd}{\sigma} + \Oh(1)\Big) \ = \ \infty. $\\
Note that actually  the difference of the formula in the last line compared to  \eqref{eq:smooth}, with the dispersion parameter $\sigma$ replacing $\sin \theta$, is negligible.
\item[(SA2)]  $\displaystyle \lim_{\sigma\to \infty} H(N,d,\sigma)=\ln(Nd)+2(\ln 2+1)$,
and we recover the well-known, average-case analysis, bound for
$\Exp_{x\in\IS^N}\ln(\msC(x))$
(see~\cite{Demmel88} and~\cite[Theorem~21.1]{Condition}).
\end{description}
\eproof

\subsection{Local analysis}

The main result of this section is the following.
\begin{theorem}\label{thm:main-local}
Let $\msC$ be a conic condition number on $\R^{N+1}$ with $N\ge 6$, with set of
ill-posed
inputs $\Sigma$. Assume that $\Sigma$ is contained in a real algebraic
hypersurface, given as the zero set of a homogeneous polynomial of
degree $d$. Let $\bar x\in \IS^N$
 and $\sigma\geq 0$. Then,
 there is an  explicit  bound $H(N,d,\sigma, \msC(\bar x))$ --defined in \eqref{eq:GaussianBound} below--  such that
$$ \Exp_{x\sim N(\bar{x}, \sigma^2\Id)}\ln \msC(x)\le H(N,d,\sigma,\msC(\bar x) ).$$
This bound satisfies   {\bf (LA0)}, {\bf (LA1)} and {\bf (LA2)}.
\end{theorem}

In order to prove Theorem~\ref{thm:main-local} we need  the following lemma.
\begin{lemma}\label{lem:otra-cota}
Assume $N\geq 2$. For all $t\in[0,\pi/2]$,
$$
    \int_t^{\frac{\pi}{2}}f(\theta)\diff\theta \leq
    \min\Big\{1,\frac{2\pi\sigma\sqrt {N+1}}{(1-e^{-\frac{1}{2\sigma^2}})t}\Big\}.
$$
\end{lemma}

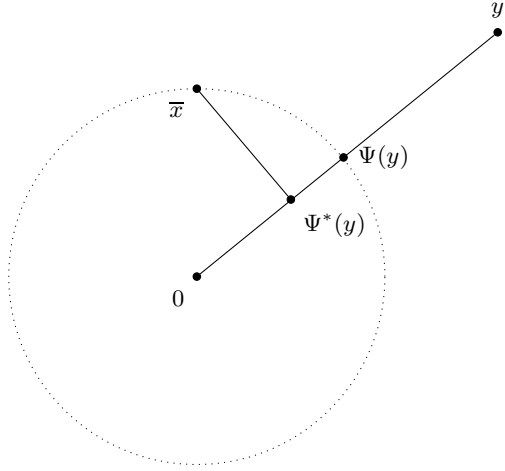
\begin{figure}[H]\centering
\begin{tikzpicture}[scale=.5,point/.style={draw,minimum size=0pt,
	inner sep=1pt,circle,fill=black}]
\draw (0,0) node(x) [point,label=220:{\footnotesize$\bar{x}$}] {};
\draw [dotted] (0,-5) circle [radius=5];
\draw (0,-5) node(q) [point,label=-120:{\footnotesize$0$}] {};
\draw (8,1.5) node(p) [point,label=-270:{\footnotesize$y$}] {};
\draw (q) -- (p);
\draw (2.5,-2.95) node(r) [point,label=-60:{\footnotesize$\Psi^*(y)$}] {};
\draw (0,0) --(r);
\draw (3.9,-1.83) node(s) [point,label=0:{\footnotesize$\Psi(y)$}] {};
\end{tikzpicture}
\caption{{\small The point $\Psi^*(y)$.}}
\label{fig:proj}
\end{figure}

\proof
The idea is to apply Markov's inequality~(e.g. \cite[Corollary~2.9]{Condition})
to the density $f$
to deduce that
$$
\int_t^{\frac{\pi}{2}} f(\theta)\diff \theta=\Prob_{\theta \sim f}(\theta \ge t)
\le \frac{1}{t}\Exp_{\theta \sim f} (\theta)
$$
Therefore we need  to bound $\displaystyle{\Exp_{\theta\sim f}(\theta)}$.
We first prove that
\begin{equation}\label{eq:des1}\Exp_{\theta\sim f}(\theta)
  \le \frac{\sqrt 2 \pi}{1-e^{\frac{-1}{2\sigma^2}}}\Exp_{y\in N(\bar{x},\sigma^2\Id)}(\|\Psi(y)-\bar{x}\|),
\end{equation}
where $\Psi$ is given by~\eqref{eq:Psi}, and then that
\begin{equation}\label{eq:des2} \Exp_{y\in N(\bar{x},\sigma^2\Id)}(\|\Psi(y)-\bar{x}\|)
\le {\sqrt{2}}\, \sigma\sqrt{N+1}.
\end{equation}
This implies
$$
\Exp_{\theta\sim f}(\theta)\le \frac{2\pi \sigma \,\sqrt{N+1}}{1-e^{\frac{-1}{2\sigma^2}}}.
$$
To show \eqref{eq:des1} we apply Proposition~\ref{prop:tech} with
$F(y)=\|\Psi(y)-\bar{x}\|$ and get
\begin{equation}\label{eq:exp1}
 \Exp_{y\in N(\bar{x},\sigma^2\Id)}(\|\Psi(y)-\bar{x}\|)
 \geq  (1-e^{\frac{-1}{2\sigma^2}})
  \Exp_{\theta \sim f}\Big(\Exp_{x\in B_{\IS}(\bar{x},\theta)}(\|x-\bar{x}\|)\Big).
\end{equation}
We claim that
\begin{equation}\label{eq:exp3}
   \Exp_{x\in B_{\IS}(\bar{x},\theta)}(\|x-\bar{x}\|) \geq \frac{\sqrt 2}{2\pi}\theta.
\end{equation}
\begin{comment}
\begin{eqnarray*}
  \Exp_{x\in B_{\IS}(\bar{x},\theta)}(\|x-\bar{x}\|) &\geq &
  \Exp_{x\in B_{\IS}(\bar{x},\theta)} \sphericalangle(x,\bar{x})\\
  &=& \frac{1}{\sfv(\theta)}
  \bigg(\int_{B_{\IS}(\bar{x},\frac{\theta}{2})}\sphericalangle[x,\bar{x}] \diff x +
  \int_{B_{\IS}(\bar{x},\theta)\setminus B_{\IS}(\bar{x},\frac{\theta}{2})}
  \sphericalangle[x,\bar{x}] \diff x \bigg)\\
  &\ge& \frac{1}{\sfv(\theta)}
  \int_{B_{\IS}(\bar{x},\theta)\setminus B_{\IS}(\bar{x},\frac{\theta}{2})}
  \frac{\theta}{2} \diff x \;=\;
  \frac{\sfv(\theta)-\sfv(\frac{\theta}{2})}{\sfv(\theta)}\frac{\theta}{2}.
\end{eqnarray*}
\end{comment}
Indeed, for  $0\le \alpha:=\sphericalangle(x,\bar x) \le \frac{\pi}{2}$, one
has
$$
   \frac{2\sqrt 2}{\pi} \alpha \le \|x-\bar x\|\le \alpha.
$$
\begin{comment}
the right-hand side inequality is quite obvious (and well-known) since
$$\|x-\bar x\|^2 = 2 (1-\cos \alpha) = 4 \sin^2 \frac{\alpha}{2},$$ which holds since  $\sin \alpha\le \alpha$ for $0\le \alpha \le \frac{\pi}{4}$.
We obtain the left-hand side inequality by showing that  $2\sin \frac{\alpha}{2}\ge \frac{2\sqrt 2}{\pi}\alpha$, i.e.
$$ \sin \alpha - \frac{2\sqrt 2}{\pi}\alpha \ge 0 \ \mbox{ for } \ 0\le \alpha\le \frac{\pi}{4}.$$
We note that in  $0$ and in $\frac{\pi}{4}$ this function gives  $0$ and furthermore, it is concave, which shows the left-hand side inequality.
\end{comment}
Therefore, writing $\sfv(\theta):=\vol(B_{\IS}(\bar{x},\theta))$,
\begin{eqnarray*}
  \Exp_{x\in B_{\IS}(\bar{x},\theta)}(\|x-\bar{x}\|) &{\geq }&
  \frac{2\sqrt 2}{\pi} \Exp_{x\in B_{\IS}(\bar{x},\theta)} \sphericalangle(x,\bar{x})\\
  &=& \frac{2\sqrt 2}{\pi\,\sfv(\theta)}
  \bigg(\int_{B_{\IS}(\bar{x},\frac{\theta}{2})}\sphericalangle[x,\bar{x}] \diff x +
  \int_{B_{\IS}(\bar{x},\theta)\setminus B_{\IS}(\bar{x},\frac{\theta}{2})}
  \sphericalangle[x,\bar{x}] \diff x \bigg)\\
  &\ge& \frac{2\sqrt 2}{\pi\sfv(\theta)}
  \int_{B_{\IS}(\bar{x},\theta)\setminus B_{\IS}(\bar{x},\frac{\theta}{2})}
  \frac{\theta}{2} \diff x \;=\;
 \frac{\theta\sqrt 2}{\pi}\left( \frac{\sfv(\theta)-\sfv(\frac{\theta}{2})}{\sfv(\theta)}\right).
\end{eqnarray*}
Now, for  $0\le \theta \le \frac{\pi}{2}$, we have $$\sin\theta =2\sin\frac{\theta}{2}\cos\frac{\theta}{2}\ge \sqrt 2 \sin\frac{\theta}{2},$$ which  implies
$$\sin \frac{\theta}{2}\le \frac{\sin \theta}{\sqrt 2}.$$
Using~\eqref{eq:volBS} twice we have, for $N\geq 6$,
$$
\sfv\Big(\frac{\theta}{2}\Big) \leq \frac{\Oh_N}{2}\Big(\sin\frac{\theta}{2}\Big)^N
\le \frac{\Oh_N}{2}\frac{1}{2^{\frac{N}{2}}}(\sin\theta)^N
\le \frac{\Oh_N}{2} \frac{1}{\sqrt{2\pi ( N+1)}}(\sin\theta)^N
\le \frac{\sfv(\theta)}{2}
$$
and we deduce that
$\frac{\sfv(\theta)-\sfv(\frac{\theta}{2})}{\sfv(\theta)}\ge\frac12$.
With this,
$$
   \Exp_{x\in B_{\IS}(\bar{x},\theta)}(\|x-\bar{x}\|) \geq \frac{\sqrt 2}{2\pi}\theta
$$
which shows \eqref{eq:exp3}.
From~\eqref{eq:exp1} and~\eqref{eq:exp3} it follows that
$$
    \Exp_{y\in N(\bar{x},\sigma^2\Id)}(\|\Psi(y)-\bar{x}\|)
    \geq \frac{\sqrt 2 (1-e^{\frac{-1}{2\sigma^2}})}{2\pi} \Exp_{\theta\sim f}(\theta),
$$ which shows \eqref{eq:des1}. We now show \eqref{eq:des2}.
We let $\Psi^*(y)$  be the closest point to
$\bar{x}$ on the line through $0$ and $y$ (see Figure~\ref{fig:proj})
and  have
\begin{eqnarray*}
\Exp_{y\in N(\bar{x},\sigma^2\Id)}(\|\Psi(y)-\bar{x}\|)
&\leq& {\sqrt{2}} \Exp_{y\in N(\bar{x},\sigma^2\Id)}(\|\Psi^*(y)-\bar{x}\|) \\
&\leq& {\sqrt{2}} \Exp_{y\in N(\bar{x},\sigma^2\Id)}(\|y-\bar{x}\|)
\;\le \; {\sqrt{2}}\, \sigma\sqrt{N+1},
\end{eqnarray*}
where the last inequality is a consequence
of~\cite[Prop. 2.10 \& Lem. 2.15]{Condition}.

This shows \eqref{eq:des2}. Therefore,
 $$
\Exp_{\theta\sim f}(\theta)\le \frac{2\pi \sigma \,\sqrt{N+1}}{1-e^{\frac{-1}{2\sigma^2}}}.
$$ as desired,
and hence,
\begin{equation}\tag*{\qed}
   \int_t^{\frac{\pi}{2}}f(\theta)\diff\theta \le \frac{2\pi\sigma\sqrt {N+1}}{(1-e^{\frac{-1}{2\sigma^2}})t}.
\end{equation}

\proofof{Theorem~\ref{thm:main-local}}
Let $t:=\arcsin\frac{1}{2\msC(\bar{x})}$. Since $\msC(\bar{x})\ge 1$,
$\frac{1}{2\msC(\bar{x})}\le \frac12$ and we have
$t\le \frac{\pi}{6}$. For all $\theta \leq t$ and
all $x\in B_{\IS}(\bar{x},\theta)$ we have
$$
\frac{1}{\msC(x)}=d_{\sin}(x,\Sigma)\ge d_{\sin}(\bar{x},\Sigma) -d_{\sin}(x,\bar{x})
\ge \frac{1}{\msC(\bar{x})}-\sin\theta \geq \frac{1}{\msC(\bar{x})}-\frac{1}{2\msC(\bar{x})}
\geq \frac{1}{2\msC(\bar{x})},
$$
which implies $\ln (\msC(x) ) \le \ln ( 2\msC(\bar{x}) )$.\\
We apply
Proposition~\ref{prop:tech} to $F(y)=\ln\msC(y)$ and use the previous  inequality and the
bounds~\eqref{eq:smooth} and~\eqref{eq:smooth-pi} to obtain
\begin{align}\label{eq:final}
  \Exp_{y\sim N(\bar{x},\sigma)} \ln\msC(y)
 &= (1-e^{\frac{-1}{2\sigma^2}})\Exp_{\theta \sim f}\Big(\Exp_{x\in B_{\IS}(\bar{x},\theta)}\ln\msC(x)\Big) +  e^{-\frac{1}{2\sigma^2}} \Exp_{x\in \IS_+^N(\bar x)} \ln\msC(x)\nonumber\\
&\le (1-e^{\frac{-1}{2\sigma^2}})\Big( \ln(2\msC(\bar{x}))\int_0^t f(\theta)\diff \theta +
  \int_t^{\frac{\pi}{2}}\Big(\ln\Big(\frac{Nd}{\sin\theta}\Big)+K\Big) f(\theta)
  \diff \theta\Big)\nonumber\\
  &\quad + e^{-\frac{1}{2\sigma^2}}(\ln(Nd)+K)\nonumber\\
 & \le  \ln \msC(\bar{x})(1-e^{\frac{-1}{2\sigma^2}})\int_0^t f(\theta)\diff \theta
    + (1-e^{\frac{-1}{2\sigma^2}})\int_t^{\frac{\pi}{2}}\ln\Big(\frac{1}{\sin\theta}\Big) f(\theta)
    \diff \theta\\
    &\;+ \ln(Nd) \Big(e^{-\frac{1}{2\sigma^2}}+(1-e^{\frac{-1}{2\sigma^2}})\int_t^{\frac{\pi}{2}} f(\theta)
    \diff\theta\Big) + K,\nonumber
\end{align}
since $$(1-e^{\frac{-1}{2\sigma^2}})\big(\ln 2 \int_0^t f(\theta)\diff\theta  + K \int_t^{\frac{\pi}{2}} f(\theta )\diff\theta \big)+ Ke^{-\frac{1}{2\sigma^2}} \le K.$$
We next bound each of the first three terms in the right-hand side.

Applying Lemma~\ref{lem:B} and the inequality $\sin(2t)\leq 2\sin t$
we obtain
\begin{eqnarray*}
(1-e^{\frac{-1}{2\sigma^2}})\int_0^t f(\theta)\diff \theta
&\le& \min\Big\{1- e^{-\frac{1}{2\sigma^2}},\frac12 (\sin(2t))^N +
  \frac{(\sin t)^N}{\sigma^{N+1}}\Big\}\nonumber \\
&\le& \min\Big\{1- e^{-\frac{1}{2\sigma^2}},\frac1{2(\msC(\bar{x}))^N} +
\frac{1}{(2\msC(\bar{x}))^N\sigma^{N+1}}\Big\}\\
&=& \min\Big\{1- e^{-\frac{1}{2\sigma^2}},\frac1{(2\msC(\bar{x}))^N} \Big(
2^{N-1}+\frac{1}{\sigma^{N+1}}\Big)\Big\}
.
\end{eqnarray*}
This bounds the first term in~\eqref{eq:final} by
\begin{equation}\label{eq:term1}
\ln \msC(\bar{x})  \min\Big\{1- e^{-\frac{1}{2\sigma^2}},\frac1{(2\msC(\bar{x}))^N} \Big(
2^{N-1}+\frac{1}{\sigma^{N+1}}\Big)\Big\}.
\end{equation}

Second, by Lemma~\ref{lem:A2}  since $t\le \frac{\pi}{6}$,
Lemma~\ref{lem:otra-cota} and $t\ge \sin t=\frac{1}{2\msC(\bar{x})}$,
\begin{align}\label{eq:2do}
\int_t^{\frac{\pi}{2}}\ln&\Big(\frac{1}{\sin\theta}\Big) f(\theta)
\diff \theta
\le \ln  \sqrt 2
+ \int_{\frac{1}{2 \msC(\bar x)}}^{\frac{\sqrt 2}{2}}  \Big(\int_{t}^{\arcsin s}
 f(\theta)\diff \theta\Big)\frac1s\diff s \nonumber\\
 \le\; & \ln  \sqrt 2
 + \int_{\frac{1}{2 \msC(\bar x)}}^{\frac{\sqrt 2}{2}}
 \min\Big\{1-e^{-\frac1{2\sigma^2}},\frac{2\pi\sigma\sqrt{N+1}}{t}\Big\} \frac{1}{s}\diff s
 \nonumber \\
  \le\; & \ln  \sqrt 2
 + \min\Big\{1,\frac{2\pi\sigma\sqrt{N+1}}{(1-e^{-\frac1{2\sigma^2}})t}\Big\}
 \Big(\ln \Big(\frac{\sqrt 2}{2}\Big) - \ln \frac{1}{ 2 \msC(\bar x)}\Big)\nonumber
 \\
  =\; & \ln  \sqrt 2
 + \min\Big\{1,\frac{4\pi\sigma\msC(\bar x)\sqrt{N+1}}{1-e^{-\frac1{2\sigma^2}}}\Big\}
\ln(\sqrt 2 \msC(\bar x)  )\nonumber
 \\
 \le\; &\min\big\{ 1,
  \frac{4\pi \sigma \msC(\bar x) \sqrt {N+1}}{1-e^{-\frac1{2\sigma^2}}}\big\}\ln \msC(\bar x) + \ln 2.
\end{align}
Also, as $t\geq 0$, we have by Lemma~\ref{lem:int1} that
\begin{eqnarray*}
 \int_t^{\frac{\pi}{2}}\ln\Big(\frac1{\sin\theta}\Big)f(\theta)
 \diff \theta
 &\le&
 \frac{1}{N}
 \Big( 1 +  \ln\big( 2^{N-1}+\frac {1}{\sigma^{N+1}}\big)-
 \ln( 1-e^{-\frac{1}{2\sigma^2}})     \Big)\\
 &\le & \frac{1}{N}
 \Big(\ln\big( 2^{N-1}+\frac {1}{\sigma^{N+1}}\big)-
 \ln( 1-e^{-\frac{1}{2\sigma^2}})\Big)   +\ln 2.
\end{eqnarray*}
Putting together this inequality and~\eqref{eq:2do} we deduce that
the second term in~\eqref{eq:final} is bounded by
\begin{align}\label{eq:term2}
  \min&\bigg\{(1-e^{-\frac{1}{2\sigma^2}})\ln(\msC(\bar x)),
 4\pi \sigma \msC(\bar x) \sqrt {N+1}\ln(\msC(\bar x)), \nonumber\\
  &\qquad
  \frac{(1-e^{-\frac{1}{2\sigma^2}})}{N}\Big(\ln\big( 2^{N-1}+\frac {1}{\sigma^{N+1}}\big)-
 \ln( 1-e^{-\frac{1}{2\sigma^2}})\Big)\bigg\}+\ln 2.
\end{align}
Finally, using again Lemma~\ref{lem:otra-cota} and
$t\ge \sin t=\frac{1}{2\msC(\bar{x})}$
we obtain
\begin{eqnarray*}\label{eq:3ro}
e^{-\frac{1}{2\sigma^2}}+( 1-e^{-\frac{1}{2\sigma^2}})\int_t^{\frac{\pi}{2}} f(\theta)\diff\theta
&\leq& e^{-\frac1{2\sigma^2}}+
\min\Big\{1-e^{-\frac1{2\sigma^2}},\frac{2\pi\sigma\sqrt{N+1}}{t}\Big\}\nonumber \\
&\leq &
\min\Big\{1,e^{-\frac1{2\sigma^2}}+4\pi\msC(\bar{x})\sigma\sqrt{N+1}\Big\}
\end{eqnarray*}
which bounds the third term in~\eqref{eq:final} by
\begin{equation}\label{eq:term3}
  \ln(Nd)\min\Big\{1,e^{-\frac1{2\sigma^2}}+4\pi\msC(\bar{x})\sigma\sqrt{N+1}\Big\}.
\end{equation}
Combining~\eqref{eq:term1},~\eqref{eq:term2}
and~\eqref{eq:term3} with the bound in~\eqref{eq:final}, we obtain
 \begin{align}\label{eq:GaussianBound}
 H(N,d,\sigma,\msC(\bar x))
 &= \ln \msC(\bar{x})  \min\Big\{1- e^{-\frac{1}{2\sigma^2}},\frac1{(2\msC(\bar{x}))^N} \Big(
2^{N-1}+\frac{1}{\sigma^{N+1}}\Big)\Big\}\nonumber \\
 & +\; \min\bigg\{(1-e^{-\frac{1}{2\sigma^2}})\ln(\msC(\bar x)),
   4\pi \sigma \msC(\bar x) \sqrt {N+1}\ln \msC(\bar x),\\
 &\qquad   \frac{(1-e^{-\frac{1}{2\sigma^2}})}{N}\Big(\ln\big( 2^{N-1}+\frac {1}{\sigma^{N+1}}\big)-
 \ln( 1-e^{-\frac{1}{2\sigma^2}})\Big)\bigg\}\nonumber \\
 &+\; \ln(Nd)\min\Big\{1,e^{-\frac1{2\sigma^2}}+4\pi\msC(\bar{x})\sigma\sqrt{N+1}\Big\}
 + \bar{K}, \nonumber
\end{align}
where $\bar{K}=\ln 2 + K=3\ln 2+2$.
We now verify that $ H(N,d,\sigma,\msC(\bar x))$ satisfies {\bf (LA0)}, {\bf (LA1)} and {\bf (LA2)}.
\begin{description}
\item[(LA0)]When
$\msC(\bar{x})=\infty$ we get
$$
H_\infty(N,d,\sigma)= \frac{(1-e^{-\frac{1}{2\sigma^2}})}{N}\Big(\ln\big( 2^{N-1}+\frac {1}{\sigma^{N+1}}\big)-
 \ln( 1-e^{-\frac{1}{2\sigma^2}})\Big)
 + \ln(Nd) + \Oh(1),
$$ which is that of \eqref{SmoothGaussianBound} (with a slightly
bigger constant) as required in {\bf (LA0)}.
\item[(LA1)]  When $\sigma\to0$, we have
$$\lim_{\sigma\to 0} H(N,d,\sigma,\msC(\bar x))=\ln(\msC(\bar{x}))+ \bar{K},$$
as required.
\item[(LA2)] Also, when $\sigma\to \infty$, we get
$$\lim_{\sigma\to \infty} H(N,d,\sigma,\msC(\bar x))= \ln(Nd)+\bar{K},$$
and we recover the average-case analysis  bound for
$\Exp_{x\in\IS^N}\ln(\msC(x))$.
\end{description}
\eproof

\end{document}